\documentclass[12pt]{amsart}

\usepackage{amsmath,amssymb,amsthm,amscd}
\numberwithin{equation}{section} 

\def\R{{\mathbb R}}

\def\cA{{\mathcal A}}

\def\cB{{\mathcal B}}

\def\cF{{\mathcal F}}

\def\cK{{\mathcal K}}

\def\cX{{\mathcal X}}

\newtheorem{prop}{Proposition}[section]
\newtheorem{theo}[prop]{Theorem}
\newtheorem{lemma}[prop]{Lemma}
\newtheorem{cor}[prop]{Corollary}
\newtheorem{rmk}[prop]{Remark}

\newtheorem{defi}[prop]{Definition}

\newtheorem{q}[prop]{Question}

\def\R{\mathbb{R}}

\let\lra=\longrightarrow

\def\mapright\#1{\,\smash{\mathop{\lra}\limits^{\#1}}\,}

\begin{document}
\bibliographystyle{plain}
\title{The Calabi flow on toric Fano Surfaces}
\date{}
\author{Xiuxiong Chen}
\author{Weiyong He}

\thanks{The first named author is partially supported by NSF. The second named author is partially supported by the PIMS postdoc fellowship.}
\maketitle

\section{Introduction}
This is a continuation of the earlier work  by the authors on the
Calabi flow \cite{chen-he, chen-he-1}. We follow the
setup of \cite{chen-he-1}, in particular we shall use the result of  formation of singularities along the Calabi flow on K\"ahler surfaces in \cite{chen-he-1}.   The readers are encouraged to
consult \cite{chen-he-1} for the setup and for the references on
this topic. The search of extremal K\"ahler metrics is now a very hot
topic in K\"ahler geometry and many people have been
contributing in this effort; we just list only a few references  \cite{Calabi82, LeBrun-Simanca, lebrun-singer, AP1, apocalgau01, APS, Fine} etc. 

We believe that the Calabi flow is a very effective tool to approach the existence of extremal metrics on compact K\"ahler manifolds.   One of the main problems on the Calabi flow is the longtime existence. 
In \cite{chen-he}, we proved that the Calabi flow exists as long as Ricci curvature tensors of the evolve metrics stay bounded.  This
is the first attempt  to understand a conjecture by the first
named author:  starting from any smooth K\"ahler metric on a
compact K\"ahler manifold (complex dimension of $n\geq 2$), the Calabi flow exists for all positive time.   In
\cite{chen-he-1}, we focused on the study of  the Calabi flow  on K\"ahler surfaces with the assumption that the Sobolev constants of the evolved metrics are uniformly bounded.  
First we  \cite{chen-he-1} studied the formation of singularity on K\"ahler surfaces.  If the curvature tensor blows up along the Calabi flow,  we could then construct a singular model, called a  {\it maximal bubble}, which is a complete asymptotically locally Euclidean (ALE) scalar flat K\"ahler surface.
Then we studied some examples where such a bubble cannot be formed;
in particular, we considered a  family of K\"ahler classes on K\"ahler surfaces of differential type of  $\mathbb{CP}^{2} \sharp
k \overline{\mathbb{CP}^{2}} (1\leq k \leq 3)$.  These surfaces are known as del Pezzo surfaces with toric
symmetry.  We  then followed the approach in
\cite{chen-lebrun-weber} to analyze all possible maximal
bubbles. Actually a maximal bubble can only  be formed in a fairly restricted way, in particular with the toric symmetry. With the aid of special geometry of manifolds we considered, in particular the toric symmetry and the discrete symmetry that those K\"ahler classes admit, we could rule out the formation of maximal bubble.  Hence we \cite{chen-he-1} could prove the longtime existence and convergence of the Calabi flow for those examples. 
However,  the analysis there is quite delicate, complicated and sometime outright challenging. It is also very hard to push these ideas
beyond the examples we considered in \cite{chen-he-1}, for example for the K\"ahler classes without discrete symmetry.

In this note we shall adopt a different strategy to rule out possible bubbles; in particular  we shall use the toric condition in a more essential way. This allows us to prove 
some longtime existence and convergence results in a fairly large family of K\"ahler classes on toric Fano surfaces. 

Let $(M, J)$ be a compact K\"ahler surface and let $[\omega]$ be a fixed K\"ahler class on $M$. We shall  use $c_1$ to denote the first Chern class of $(M, J)$.  We may define a functional
\[
\cB([\omega])=32\pi^2\left(c_1^2+\frac{1}{3}\frac{(c_1\cdot[\omega])^2}
{[\omega]^2}\right) +\frac{1}{3}\|\cF\|^2,
\]
where  $\|\cF\|^2 $ is the
norm of Calabi-Futaki invariant \cite{Futaki-Mabuchi, chen05}. Our main result is  
\begin{theo}\label{T-M} Let $(M, [\omega], J)$ be a toric Fano
surface with positive extremal Hamiltonian potential. If the Calabi flow
initiates from a K\"ahler metric with toric symmetry satisfying
\begin{equation}\label{E-1-1}
\int_MR^2dg<\cB([\omega]),
\end{equation}
then the Calabi flow exists for all time and converges subsequently to an
extremal metric in $[\omega]$ in Cheeger-Gromov sense.
\end{theo}
The definition of {\it extremal Hamiltonian potential} will be given in Section
2. An immediate corollary of Theorem \ref{T-M} is:

\begin{cor}
Let  $(M, [\omega], J)$ be a toric  Fano surface  with
positive extremal Hamiltonian potential.  If there is  a toric metric $\omega_0\in [\omega]$ such that the Calabi energy of $\omega_0$
is less than $\cB([\omega])$, 
then there exists an extremal metric in $[\omega]$.
\end{cor}

\begin{rmk}  In  \cite{Dona05, Dona06},  Donaldson used
a continuity method to deform metrics to seek extremal metrics on toric surfaces and has
made striking progress on existence of constant scalar curvature  metrics. His
approach uses convex analysis, which depends on the fact that, one can express a toric metric on a toric surface  in terms of  a convex function in
a convex polytope in $\mathbb{R}^2$.
\end{rmk}

\noindent {\bf Acknowledgment:} The authors have been
lecturing on the results since the fall of 2007.  The
authors want to thank B. Wang for sharing his insights on geometric flows arising from K\"ahler geometry. The first named author wants to thank Professor S.K. Donaldson for insightful comments. The second named author wants to thank Professor A. Futaki for helpful comments on extremal Hamiltonian potential.
The present version is  rewritten thoroughly and the authors are grateful to the referee for pointing out a numerical error in a previous version and numerous  suggestions which help improve the presentation of the work.

\section{Sobolev Constant}
In this section we shall prove that the Sobolev constants of the evolved metrics along the Calabi flow on Fano surfaces are uniformly bounded under certain natural geometric conditions. 
We shall first define  an  extremal Hamiltonian potential of an {\it invariant  K\"aher metric} in a
fixed K\"ahler class $(M, [\omega])$, which is essentially  given by \cite{Futaki-Mabuchi, Mabuchi}.
Recall that an extremal vector
field for $(M, [\omega])$ is {\it a priori} determined \cite{Futaki-Mabuchi} up to
conjugation. Let $\cX$ be the extremal vector field and let $\cX_\R$ be the real part of $\cX$.
Define $\cK_\cX$ to be the set of  all the invariant metrics in $[\omega]$ which satisfy
\[
\cK_\cX=\{\omega: L_{\cX_\R} \omega=0\}.
\]
For any
K\"ahler metric $\omega\in \cK_\cX$, one can  define the real potential $\theta_{ \omega}$ \cite{Mabuchi} by
\[
\nabla_\omega \theta_\omega=\cX_\R
\]
which satisfies the normalized condition
\[\int_M\theta_{ \omega}\omega^n=0.
\] 
We have the $L^2$
orthogonal decomposition \cite{hwang, chen05}
\[
R_\omega=\underline{R}+\theta_{
\omega}+\theta_{\omega}^{\bot},
\]
where the average of the scalar curvature $\underline{R}$ is determined by $(M, [\omega])$. We can then define the extremal Hamiltonian potential as
\begin{defi}
Let $\omega\in \cK_\cX$, the  extremal Hamiltonian potential of $\omega$  is given by 
\[
\rho_\omega=\underline{R}+\theta_{\omega}.
\]
\end{defi}

An extremal metric $\omega$ then satisfies 
$R_\omega=\rho_\omega$. 
By its definition  \cite{Mabuchi}, the maximum and the minimum of
$\theta_{ \omega}$ are the invariants of $(M, [\omega])$.  We may denote \[
\theta_{-}=\min_{\omega\in \cK_\cX} \theta_{\omega}, \;
\theta_{+}=\max_{\omega\in \cK_\cX} \theta_{ \omega}.
\]
We can also denote \[\rho_{-}=\underline{R}+\theta_{-},\; \rho_{+}=\underline{R}+\theta_{+}.\]
It is clear that $\rho_{-}$ and $\rho_{+}$ are the minimum and maximum of the scalar curvature of an extremal metric  respectively if it exists in $[\omega]$.
If $(M, [\omega])$ is a Fano surface and  the Futaki invariant  of $[\omega]$ is zero, then
$\theta_{-}=\theta_{+}=0$ and so $\rho_{-}=\rho_{+}=\underline{R}$  is positive. Hence $\rho_{-}$ is positive for Fano surfaces of differential type of $\mathbb{CP}^2\sharp
k\overline{\mathbb{CP}^2} (4\leq k\leq 8)$.  
An interesting question is
\begin{q}\label{P-S-1}Let $(M, [\omega])$ be a  toric Fano surface, is $\rho_{-}$ positive?
\end{q}
Note that $\rho_{-}$ is an invariant of $(M, [\omega])$ and it is computable in particular when $(M, J)$ is a toric Fano surface. Hence one can check numerically whether $\rho_{-}$ is positive or not for any given K\"ahler class on $M$. 
However, in general it seems not easy to verify that it is positive since its expression is quite complicated.
Without giving detailed argument, S. Simanca claimed  that the answer to Question \ref{P-S-1}
to be correct (cf. \cite{simanca2}). Since no detailed computation
is given by Simanca in \cite{simanca1, simanca2}, we feel it
is important  to point out  why such a
statement is plausible. Note that the average of the
scalar curvature on $(M, [\omega])$ is positive when $M$ is a  Fano surface.  Intuitively, if there is an extremal metric, then the scalar curvature  
of the extremal metric should be positive since it minimizes the Calabi energy.
To verify this, one needs to consider the case of $M\sim\mathbb{CP}^2\sharp
k\overline{\mathbb{CP}^2}$ $(k=1, 2, 3)$. 
When $k=1$, one can check the scalar curvatures of all the extremal metrics constructed by E. Calabi \cite{Calabi82} are positive.   LeBrun-Simanca \cite{LeBrun-Simanca}  computed the Futaki invariant and the extremal vector field of a K\"ahler class explicitly for K\"ahler surfaces with a semi-free $\mathbb{C}^{*}$ action. In particular their results can be applied to toric Fano surfaces and one can compute further $\rho_{-}$. For example some explicit formula is given  in \cite{zhou-zhu}. However, it seems that only when $(M, [\omega])$ admits some additional discrete symmetry, the formula of $\rho_{-}$ is simple enough and  one can check directly that it is actually  positive.
For example, when $k=2$, it is proved that $\rho_{-}$ is positive for the bilaterally symmetric K\"ahler classes 
\cite{chen-lebrun-weber}.

We shall then  show how to bound the Sobolev
constants on Fano surfaces under natural geometric conditions. The
idea dates back to Tian \cite{tian} for K\"ahler metrics of  constant scalar
curvature (see \cite{Tian-Viaclovsky} also) and it is generalized to  extremal
metrics in Chen-Weber \cite{chen-weber}. 
\begin{lemma}\label{L-2-1} Let $(M, [\omega])$ be a Fano 
surface such that $\rho_{-}>0$ and let $g$ a K\"ahler
metric in $[\omega]$. If  $g$ is invariant ($g\in \cK_\cX$) and
\begin{equation}\label{E-2-3}
\int_MR^2dg<\cB([\omega]),
\end{equation}
then the Sobolev constant is bounded a priori as in \eqref{E-2-13}. \end{lemma}

When $g$ is not invariant, Lemma \ref{L-2-1} still holds with
stronger restriction on the Calabi energy. But we shall not need this.  We define the
Sobolev constant  for a compact $4$ manifold  $(M, g)$ to be the smallest constant $C_s$ such that the
estimate holds, for any $f\in W^{2, 2}(M, g)$, 
\begin{equation}\label{E-2-4}
\|f\|^2_{L^4}\leq C_s\left(\|\nabla
f\|_{L^2}^2+V^{-1/2}\|f\|^2_{L^2}\right),
\end{equation}
where $V$ is the volume of the manifold  $(M, g)$. Note the Sobolev
inequality (\ref{E-2-4}) is scaling-invariant. When the Yamabe
constant is positive, the Sobolev constant is essentially bounded
by the Yamabe constant \cite{aubin}. Recall that the Yamabe
constant for a conformal class $[g]$ of Riemannian metrics on a
compact 4 manifold is given by
\[
Y_{[g]}=\inf _{\tilde{g}\in
[g]}{\int_MR_{\tilde{g}}d\tilde{g}\over \sqrt{\int_Md\tilde{g}}}.
\]
By the celebrated work of Trudinger, Aubin and Schoen
\cite{aubin1, LeeParker}, for any conformal class $[g]$ the
infimum is achieved by the so-called Yamabe minimizer $g_Y\in [g]$
which necessarily has constant scalar curvature. If $\tilde
g=u^2g$, the scalar curvature is given by
\[
R_{\tilde g}=u^{-3}(6\triangle_g u +R_{g}u),
\]
so the Yamabe constant takes the formula
\begin{equation}\label{2-1}
Y_{[g]}=\inf_{u\neq 0}\frac{\int_M(6|\nabla
u|^2+R_gu^2)dg}{(\int_Mu^4dg)^{1/2}}.
\end{equation}
Now we are in the position to prove Lemma \ref{L-2-1}.
\begin{proof} 
We can rewrite \eqref{E-2-3}  as 
\begin{equation}\label{2-3}
96\pi^2c_1^2-2\int_MR^2dg> \int_M(R-\underline R)^2dg-\|\cF\|^2.
\end{equation}
Following  computation in \cite{tian, chen-lebrun-weber} (for example, see Section 5 \cite{chen-lebrun-weber}), we have
\begin{equation}\label{E-2-7}
Y^2_{[g]}\geq 96\pi^2c_1^2-2\int_MR^2dg.
\end{equation}
It then follows from \eqref{2-3} and \eqref{E-2-7} that
\begin{equation}\label{E-2-8}
Y_{[g]}^2>\int_M(R-\underline R)^2dg-\|\cF\|^2.
\end{equation}
We shall  need a decomposition formula of the Calabi energy \cite{hwang, chen05},
\begin{equation}\label{E-2-2}
\int_M(R-\underline{R})^2dg-\|\cF\|^2=\int_M(R-\underline{R}-\theta_{
\omega})^2dg.
\end{equation}
First we show that $Y_{[g]}$ has to be positive. Pick up a sequence of functions $u_i$ ($u_i\ne 0$) which minimizes the expression in \eqref{2-1}. Hence we have
\begin{equation}\label{2-2}
Y_{[g]}+\epsilon_i= \frac{\int_M(6|\nabla
u_i|^2+R_gu^2_i)dg}{(\int_Mu^4_idg)^{1/2}},
\end{equation}
such that $\epsilon_i\rightarrow 0$ when $i \rightarrow \infty$. 
We can rewrite \eqref{2-2} as
\begin{equation}\label{2-4}
(Y_{[g]}+\epsilon_i)\|u_i\|^2_{L^4}=6\int_M|\nabla
u_i|^2dg+\int_MRu_i^2dg,
\end{equation}
where we write $R=R_g$ for simplicity. 
It then follows from \eqref{2-4} that
\begin{equation}\label{E-2-9}
(Y_{[g]}+\epsilon_i)\|u_i\|^2_{L^4}-\int_M(R-\underline{R}-\theta_{
\omega})u_i^2dg=6\|\nabla u_i\|^2_{L^2}+(\underline{R}+\theta_{
\omega})\|u_{i}\|^2_{L^2}.
\end{equation}
By Cauchy-Schwarz inequality, we compute 
\begin{equation}\label{2-5}
\left|\int_M(R-\underline{R}-\theta_{\omega})u_i^2dg\right|\leq \left(\int_M (R-\underline{R}-\theta_{\omega})^2 dg\right)^{1/2}\left(\int_M u_i^4dg\right)^{1/2}.
\end{equation}
Then we compute, by \eqref{2-5}, 
\begin{equation}\label{2-6}
(Y_{[g]}+\epsilon_i)\|u_i\|^2_{L^4}-\int_M(R-\underline{R}-\theta_{
\omega})u_i^2dg\leq (Y_{[g]}+\epsilon_i+\|R-\underline{R}-\theta_{\omega}\|_{L^2})\|u_i\|_{L^4}^2. 
\end{equation}
If $Y_{[g]}<0$, then by \eqref{E-2-8} and \eqref{E-2-2}, we know that
\begin{equation}\label{2-7}
Y_{[g]}+\|R-\underline{R}-\theta_{\omega}\|_{L^2}<0. 
\end{equation}
Since $g$ is fixed,  then by \eqref{2-7}, 
$Y_{[g]}+\|R-\underline{R}-\theta_{\omega}\|_{L^2}+\epsilon_i$ is less than zero for sufficiently large $i$; hence by \eqref{2-6},   we can get that for $i$ large enough, 
\begin{equation}\label{2-8}(Y_{[g]}+\epsilon_i)\|u_i\|^2_{L^4}-\int_M(R-\underline{R}-\theta_{
\omega})u_i^2dg<0.\end{equation}
However $\underline{R}+\theta_{
\omega}\geq \rho_{-}>0$, the right hand side of \eqref{E-2-9} is then positive, which contradicts \eqref{2-8}. Hence $Y_{[g]}>0$; it then follows from \eqref{E-2-8} that
\begin{equation}\label{2-10}
Y_{[g]}> \|R-\underline{R}-\theta_\omega\|_{L^2}. 
\end{equation}
We can then rewrite \eqref{2-1} as, for $u>0$, 
\begin{equation}\label{E-2-10}
\|u\|_{L^4}^2\leq \frac{6}{Y_{[g]}}\|\nabla
u\|_{L^2}^2+\frac{1}{Y_{[g]}}\int_MRu^2dg.
\end{equation}
It is easy to see that (\ref{E-2-10}) holds for any $u$ since
$|\nabla |u||\leq |\nabla u|$ at $u\neq 0$.
Now we  rewrite (\ref{E-2-10}) as
\begin{equation}\label{2-9}
\|u\|_{L^4}^2-\frac{1}{Y_{[g]}}\int_M(R-\underline{R}-\theta_{
\omega})u^2dg\leq\frac{6}{Y_{[g]}}\|\nabla
u\|_{L^2}^2+\frac{1}{Y_{[g]}}\int_M(\underline{R}+\theta_{
\omega})u^2dg.
\end{equation}
Note that $\underline{R}+\theta_{ \omega}\leq \rho_{+}$.
It follows from \eqref{2-9} and Cauchy-Schwarz inequality, that
\begin{equation}\label{E-2-11}\left(1-\frac{1}{Y_{[g]}}\|R-\underline{R}-\theta_{
\omega}\|_{L^2}\right)\|u\|_{L^4}^2\leq\frac{6}{Y_{[g]}}\|\nabla
u\|_{L^2}^2+\frac{\rho_{+}}{Y_{[g]}}\|u\|_{L^2}^2.
\end{equation}
It then follows from (\ref{E-2-11}) that the Sobolev constant of  $g$ is bounded a priori. In other words, we can
get that
\begin{equation}\label{E-2-13}
C_s\leq \max\left\{
\frac{6}{Y_{[g]}-\|R-\underline{R}-\theta_{\omega}\|_{L^2}},\;
\frac{\sqrt{V}\rho_{+} }{Y_{[g]}-\|R-\underline{R}-\theta_{
\omega}\|_{L^2}}\right\}.
\end{equation}
\end{proof}

\section{Rule Out Bubbles}
In this section we shall prove Theorem \ref{T-M}.  
First let us recall the formation of singularity along the Calabi flow on K\"ahler surfaces. Let $(M, [\omega])$ be a toric Fano surface as in Theorem \ref{T-M}. 
Suppose that the Calabi flow exists on $[0, T)$, $0<T\leq \infty$ and the curvature tensor blows up when $t\rightarrow T$. 
Note that under the assumption in Theorem \ref{T-M}, the Sobolev constants of the evolved metrics are uniformly bounded by Lemma \ref{L-2-1}, since the Calabi energy is decreasing along the flow. Hence the result (Theorem 1.1, \cite{chen-he-1}) is applicable.
Since the blowing up process is required in the following argument, we shall  state the result as follows. 
\begin{prop}\label{P-1}Keep the assumption in Theorem \ref{T-M}. 
If the curvature blows up when $t\rightarrow T$, there exists a sequence of points $(x_i, t_i) \in (M, [0, T))$
where $ t_i \rightarrow T\;$ and
$ Q_i=\max_{t\leq t_i}|Rm|= |Rm(x_i, t_i)|
 \rightarrow \infty$ such that the pointed manifolds \[(M,
x_i,Q_i g (t_i+t/Q_i^2))\] converge locally smoothly to an ancient solution of the Calabi flow
 \[(M_\infty, x_\infty, g_\infty(t)),  t\in (-\infty,
0].\]  Moreover,  $g_\infty(t)\equiv g_\infty(0)$ and $g_\infty:=g_\infty(0)$ is a complete scalar flat ALE K\"ahler metric on $M_\infty$.\end{prop}
One of the key points in \cite{chen-lebrun-weber} is that $(M_\infty, g_\infty)$, as a limit of pointed manifolds $(M, g_i)$,  is toric since $g_i:=Q_ig(t_i)$ is toric. Moreover
$(M_\infty, g_\infty)$  contains holomorphic cycles. The result (Proposition 16, \cite{chen-lebrun-weber}) is only stated for $M\sim\mathbb{CP}^2\sharp 2\overline{\mathbb{CP}^2}$, but the result and the proof hold for all toric Fano surfaces without any change. We shall state the result as follows.
\begin{prop} \label{P-2}Keep the same assumption as in Theorem \ref{T-M}. Suppose the curvature tensor blows up along the Calabi flow  and let $(M_\infty, g_\infty)$ be a maximal bubble. Then $(M_\infty, g_\infty)$ is toric and $H_2(M_\infty, \mathbb{Z})$ is generated by holomorphically embedded $\mathbb{CP}^1$s in $M_\infty$.
\end{prop}
On the other hand, we show that a holomorphic cycle cannot be formed in such a blowup process. The idea is more lucid when the cohomology class $[\omega]$ is rational. 
\begin{prop}\label{P-3}
Keep the same assumption as in Theorem \ref{T-M}. Let $[\omega]\in H^{2}(M, \mathbb{Q})$. Then $(M_\infty, g_\infty)$ cannot contain a holomorphic $\mathbb{CP}^1$.
\end{prop}
\begin{proof} $(M_\infty, g_\infty)$ is the limit of pointed manifolds $(M, g_i)$. Hence there is a sequence of
compact set $K_i$, $K_i\subset K_{i+1}$, $\cup K_i=M_\infty$, and
a sequence of  diffeomorphisms $\Phi_i: K_i\rightarrow
\Phi_i(K_i)\subset M$,  \[\Phi_i^{*}(g_i)\rightarrow g_\infty,\] where the convergence is  smooth in $K_{i-1}$. 
Let $S$ be an embedded holomorphic $\mathbb{CP}^1$ in $M_\infty$.  There is a sequence of compact two spheres, which
are denoted as $S_i=\Phi_{i}(S)$ and  $S_i\subset \{M,  Q_i g(t_i)\}$.  Let
$\omega_\infty$ be the K\"ahler form of  $g_\infty$ and let $\omega_i=Q_i\omega(t_i)$ be  the K\"ahler form of $g_i$.
Since $\Phi_i^{*}g_i$ converges to $g_\infty$ smoothly, then for any fixed positive constant $\epsilon$ 
we have
\begin{equation}\label{E-3-1}
\left|\int_{S_i}\omega_i-\int_S\omega_\infty\right|=\left|\int_{S}\Phi_i^{*}\omega_i-\int_S\omega_\infty\right|<\epsilon
\end{equation}
when $i$ is sufficiently large. Hence $\int_{S_i}\omega_i$ is uniformly bounded  and  then
 \begin{equation}\label{3-2}\int_{S_i}\omega (t_i)=\frac{1}{Q_i}\int_{S_i} \omega_i \rightarrow 0.\end{equation} On the other hand, we know that 
\[
\int_{S_i}\omega(t_i)=\int_{S_i}\omega=[\omega][S_i]=a_i
\]
is a constant depending only on $[\omega], [S_i]$. Since $[\omega]\in H^{2}(M, \mathbb{Q})$, there exists some $k\in \mathbb{N}$ such that
$[k\omega]\in H^2(M, \mathbb{Z})$. It then follows that
$\int_{S_i} k\omega$ is an integer, hence $ka_i$ is an integer
for any $i$. By \eqref{3-2},  $a_i\rightarrow 0$, hence $ka_i$ has to be zero when $i$ large enough.  It then follows that
$a_i=0$ when $i$ is sufficiently large. If $a_i=0$, by (\ref{E-3-1}), it follows
that
\[
\int_S\omega_\infty=0.
\]
This contradicts  that $S$ is a holomorphic embedded
$\mathbb{CP}^1$ in $M_\infty.$
\end{proof}
When $[\omega]$ is not a rational class, 
the proof is more involved. The key is then to show that $\{[S_i]\}$ can only contain finite many homology classes, which rely on \eqref{E-3-1}, \eqref{3-2} and positivity of a K\"ahler class. 
\begin{prop}\label{P-4}Keep the same assumption as in Theorem \ref{T-M}. $(M_\infty, g_\infty)$ cannot contain a holomorphic $\mathbb{CP}^1$. 
\end{prop}
\begin{proof}
Keep the same notations as in Proposition \ref{P-3}. It is clear that we can still get \eqref{E-3-1} and \eqref{3-2} and when $i\rightarrow \infty$, 
\begin{equation}\label{E-3-2}[\omega][S_i]=a_i\rightarrow 0.
\end{equation} We show that any such sequence $\{[S_i]\}$ contains only
finite homology classes in $H_2(M,\mathbb{Z})$.
Recall that the self-intersection of $S\in H_2(M_\infty,
\mathbb{Z})$ is a negative integer \cite{chen-lebrun-weber}. Let
$[S][S]=-k,$
for some fixed integer  $k\geq 1.$ Since the self-intersection is invariant under diffeomorphism, hence for any $i$, 
\begin{equation}\label{E-3-3} [S_i ][S_i]=-k.
\end{equation}
The toric Fano surfaces are described as $\mathbb{CP}^2,
\mathbb{CP}^1\times \mathbb{CP}^1, \mathbb{CP}^2\sharp
\overline{\mathbb{CP}^2}$, $\mathbb{CP}^2\sharp 2
\overline{\mathbb{CP}^2}$ ($\mathbb{CP}^2$ blown up at two
distinct points), $\mathbb{CP}^2\sharp 3 \overline{\mathbb{CP}^2}$
($\mathbb{CP}^2$ blown up at three non-linear points). We only exhibit the
example when  $M\sim\mathbb{CP}^2\sharp 3 \overline{\mathbb{CP}^2}$, all
other examples are similar (and simpler). 
Let $H$ be a hyperplane in $\mathbb{CP}^2$.  $M$ can be obtained by blown up at three generic points on $\mathbb{CP}^2$. After blown up,  we still use $H$ to denote the corresponding hypersurface on $M$ and $E_i, i=1, 2, 3$ to denote the exceptional divisors. 
For simplicity, we use $[H], [E_i]$ to denote the homology classes and their Poincar\'e dual-the cohomology classes.
The K\"ahler classes on
$M$ can be expressed
as
\[
[\omega]_{x, y, z}=3[H]-x[E_1]-y[E]_2-z[E_3].
\]
Since $[\omega]$ is a positive class,  then $x, y, z$ have to satisfy that \begin{equation}\label{3-4}0<x, y, z;~\mbox{and}~ x+y, y+z, x+z<3.\end{equation}
We can see \eqref{3-4} as follows;  for example, $x=[E_1][\omega]_{x, y, z}>0$ and  $H-E_1-E_2$ is a holomorphic curve which has area $3-x-y$ with respect to $[\omega]_{x, y, z}$, hence $x+y<3$. 
And $H_2(M, \mathbb{Z})$ can be generated by $\{[H], [E_i], i=1, 2, 3\}$,  we can then express $[S_i]$ as
\[[S_i]=m[H]+n[E_1]+j[E_2]+l[E_3],
\]
for some integers $m, n, j, l$. We can write (\ref{E-3-2}) and (\ref{E-3-3}) 
as, when $i\rightarrow \infty$, 
\begin{equation}\label{E-3-5}3m-nx-jy-lz\rightarrow
0\end{equation} and 
\begin{equation}\label{E-3-4}
m^2-n^2-j^2-l^2=-k.
\end{equation}  We can compute, by \eqref{E-3-5},
\begin{equation}\label{E-3-6}
n^2+j^2+l^2\geq \frac{(nx+jy+lz)^2}{x^2+y^2+z^2}\rightarrow
\frac{9m^2}{x^2+y^2+z^2}.
\end{equation}
Hence, by \eqref{E-3-4} and \eqref{E-3-6},
\[
m^2+k+1=n^2+j^2+l^2+1\geq\frac{9m^2}{x^2+y^2+z^2}.
\]
But by \eqref{3-4}, it is easy to see that
\[
x^2+y^2+z^2<9.
\]
For any fixed $x, y, z$, it then follows that
\[
m^2\left(\frac{9}{x^2+y^2+z^2}-1\right)\leq k+1. 
\]
It follows that $m$ has at most finite many solutions. So there are at most finite many $m, n, j, l$ such that 
\eqref{E-3-5} and \eqref{E-3-4} are satisfied. It then
follows that the homology classes of $[S_i]$ are finite. Hence we
can find a subsequence $S_{\bar i}$ of $S_i$, such that 
$[S_{\bar i}]\in H_2(M, \mathbb{Z})$ has the same homology class for any $\bar i$. Hence 
$a_{\bar i}=[\omega][S_{\bar i}]$ is a
constant independent of $\bar i$. By \eqref{E-3-2}, $a_{\bar i}\equiv 0$. It then
follows that $[S][\omega_\infty]=0$ by \eqref{E-3-1}. This contradicts that $S$ is a holomorphic cycle in $M_\infty$. \end{proof}

\begin{rmk}Similar idea can be applied to the Calabi flow on toric surfaces, if one assumes that the Sobolev constants of the evolved metrics are uniformly bounded.
\end{rmk}

Now we shall state a convergence result for the Calabi flow. 
\begin{prop}\label{P-5} Let $(M, J)$ be a K\"ahler manifold. Suppose $(M, g(t), J), 0\leq t<\infty$ is a solution of the Calabi flow such that the Sobolev constants and the curvature tensors of the evolved metrics are uniformly bounded. 
Then for every sequence $t_i\rightarrow \infty$, there is a subsequence $t_{i_k}$ and a sequence of diffeomorphisms $\Phi_{i_k}: M\rightarrow M$ such that,
\[\Phi_{i_k}^{*}g(t_{i_k})\rightarrow g_{\infty}, {\Phi^{-1}_{i_k}}_{*}\circ J\circ  {\Phi_{i_k}}_{*}\rightarrow J_\infty,\]
under a fixed gauge, where the convergence is in $C^\infty$ topology and $(M, g_\infty, J_\infty)$ is an extremal K\"ahler manifold with complex structure $J_\infty$.  
\end{prop}
\begin{proof}By assumption both Sobolev constants and curvature tensors are bounded, then  all higher derivatives of curvature tensors are uniformly bounded, for example see Lemma 4.2 in \cite{chen-he-1}. It then follows from the standard ideas in Ricci flow (see Hamilton \cite{Hamilton951}) to get similar compactness results for the Calabi flow. For a sequence 
$t_i\rightarrow \infty$, there is a subsequence $t_{i_k}\rightarrow \infty$ such that
\[
\{M, g(t+t_{i_k}), -t_{i_k}\leq t\leq 0\}\rightarrow \{M_\infty, g_\infty(t), -\infty\leq t\leq 0\}
\]
in Cheeger-Gromov sense. The argument is well known in geometric flows  and we shall skip the details.  Let $g_\infty=g_{\infty}(0), g_{i_k}=g(t_{i_k})$. In particular, $(M, g_{i_k})\rightarrow (M_\infty, g_\infty)$. Namely, there exists a sequence of diffeomorphisms 
$\Phi_{i_k}: M\rightarrow M_\infty$ such that 
\[
\Phi_{i_k} ^{*}g_{i_k}\rightarrow g_\infty.\]
If necessary, by taking a subsequence, we can get that , $J_{i_k}={\Phi^{-1}_{i_k}}_{*}\circ J\circ  {\Phi_{i_k}}_{*}\rightarrow J_\infty$. 
Since $\nabla_{g_{i_k}}J_{i_k}=0$, it follows that $\nabla_{g_\infty}J_\infty=0$, hence $J_\infty$ is still a complex structure which is  compatible with $g_\infty$. 
We then show $g_\infty$ is an extremal metric. This follows from that the Calabi flow is the gradient flow of the Calabi energy. 
For any $t_0\in (-\infty, 0]$, we choose the sequence  $\{t_{i_k}\}$ such that $t_{i_{k}}<t_{i_{k+1}}+t_0$. Let  $\mathcal{C}(g)$ be the Calabi energy of $g$. Since the Calabi energy is decreasing along the Calabi flow, we have
\[
\mathcal{C}(g_\infty)=\lim_{t_{i_{k}}\rightarrow\infty} \mathcal{C}(g(t_{i_{k}}))\geq \lim_{t_{i_{k+1}}\rightarrow \infty}\mathcal{C}(g(t_0+t_{i_{k+1}}))=\mathcal{C}(g_\infty(t_0)).
\]
It then follows that $g_\infty(t)$ is an extremal metric for any $t\in (-\infty, 0]$. 
\end{proof}

\begin{rmk} In general $J_\infty$ does not have to be the same as $J$. 
\end{rmk}

Now we are in the position to prove Theorem \ref{T-M}. We argue  by contradiction.
\begin{proof}
By Lemma \ref{L-2-1},  the Sobolev constants of evolved metrics are uniformly bounded under the assumption in Theorem \ref{T-M}. If the curvature tensors are not uniformly bounded,  there is a contradiction by Proposition \ref{P-1}, \ref{P-2} and \ref{P-4}. 
Hence the curvature tensors have to be uniformly bounded and the Calabi flow exists for all time. It then follows that $(M, g(t), J)$ converges to an
extremal metric  $(M, g_\infty, J_\infty)$ subsequently  in Cheeger-Gromov sense by Proposition \ref{P-5}. We then finish the proof by showing that  $(M, J_\infty)$ is biholomorphic to $(M, J)$. The proof follows from
\cite{chen-lebrun-weber} (Theorem 27) by using the toric  condition carefully and the
classification of complex surface. Theorem 27 in \cite{chen-lebrun-weber} states only for $M\sim\mathbb{CP}^2\sharp 2\overline{\mathbb{CP}^2}$ but the proof holds for all toric Fano surfaces. The key is that in the limiting process, the torus action converges and 
 $(M, g_\infty, J_\infty)$ is still toric. Moreover, the 2-torus action for $(M,
g_\infty, J_\infty)$ is holomorphic with respect to $J_\infty$. We shall sketch the argument for  $M\sim\mathbb{CP}^2\sharp3\overline{\mathbb{CP}^2}$. The readers can refer to \cite{chen-lebrun-weber} for details.   When $M\sim\mathbb{CP}^2\sharp3\overline{\mathbb{CP}^2}$, each of holomorphic curves $H, E_1, E_2, E_3$ is the fixed point set of the isometric action of some circle action of
2-torus, and so each is totally geodesic with respect to the
metrics along the Calabi flow. By looking at the corresponding
fixed points set of the limit action of circle subgroups, we can
find corresponding  totally geodesic 2-spheres in $(M, g_\infty, J_\infty)$ which are the limits of the image of these submanifolds. Moreover, these limit 2-spheres are holomorphic with respect to $J_\infty$ and the
homological intersection numbers of these holomorphic spheres  do not vary. Namely, we have still three holomorphic $\mathbb{CP}^1$s with self-intersection $-1$ as the images of the original exceptional divisors $E_1, E_2, E_3$.  Thus, by blowing down the images of $E_1, E_2, E_3$ and applying the
classification of the complex surface, we conclude that $(M,
J_\infty)$ is  biholomorphic to $\mathbb{CP}^2$  blown up three generic points. So there exists a
diffeomorphism $\Psi$ such that $\Psi_{*} J=J_\infty$. So
$\Psi^{*}g_\infty$ is an extremal metric in the class $[\omega]$ for
$(M, J)$.
\end{proof}

\begin{rmk}We may define a functional
\[
\mathcal{A}[\omega]=\frac{(c_1\cdot[\omega])^2}{[\omega]^2}+\frac{1}{32\pi^2}\|\cF\|^2.
\]
This functional has the important property \cite{chen05, Dona01} that any K\"ahler
metric $g$ in the class $[\omega]$ satisifies the curvature inequality
\[
\int_M R^2dg\geq 32\pi^2\cA([\omega])
\]
with equality if and only if  $g$ is an extremal metric.
A necessary condition for \eqref{E-1-1} to hold is that
$(M, [\omega])$ satisfies the generalized Tian's condition
in \cite{chen-weber},
\[
c_1^2>\frac{2}{3}\cA([\omega]).
\]
\end{rmk}

\bibliographystyle{plain}

\vspace{3mm}

 Xiuxiong CHEN,  Department of Mathematics,
University of Wisconsin-Madison, Madison, WI, 53706; xxchen@math.wisc.edu.\\

 Weiyong HE, 
Department of Mathematics,
University of British Columbia, Vancouver, Canada, V6T 1Z2; whe@math.ubc.ca.
 
Current address:
Department of Mathematics, 
University of Oregon, Eugene OR, 97403;  whe@uoregon.edu.
\end{document}